\documentclass[10pt]{amsart}
\usepackage{amsmath,amsfonts,latexsym,amssymb}

\usepackage{graphicx}
\usepackage{ hyperref}
\usepackage{graphicx}

\newtheorem{theorem}{Theorem}[subsection]
\newtheorem{lemma}[theorem]{Lemma}
\newtheorem{coro}[theorem]{Corollary}
\newtheorem{prop}[theorem]{Proposition}
\newtheorem{fact}[theorem]{Fact}

\def\HH{\mathbb H}

\def\RR{\mathbb R}
\def\ZZ{\mathbb Z}

\def\MM{\Sigma}
\def\TT{\mathcal{T}(\MM)}
\def\TTE{\mathcal{T}_{<\epsilon}(\MM)}

\def\XX1{\mathcal{X}_{11}}
\def\tr{\mathrm{tr}}

\def\ei{e^{i\theta}}
\def\fntb{ \tau_\alpha^t}
\def\fnt{ \hat{\tau}_\alpha^t}

\def\BP{\partial \overline{P}}
\def\fund{\pi_1(\MM)}

\def\MM{\Sigma}
\def\HH{\mathbb{H}}
\def\dHH{\partial \mathbb{H}}

\def\axG{\mathrm{ax}(\Gamma)}
\def\ax{\mathrm{ax}}
\def\axg{\ax(\gamma)}
\def\psl{PSL_2(.)}
\def\sl{SL(2,\RR)}

\def\tw{\mathrm{tw}_\alpha}

\def\aza{\alpha_1 \angle_z \alpha_2}

\def\auta{\mathrm{Aut}(\axG)}

\def\comm{\text{Comm}}

\def\psl{\mathrm{PSL}(2,\mathbb{R})}
\def\llg{\ell_\gamma}
\def\CH{C(\Lambda)}

\def\LL{\Lambda}

\begin{document}

%
%

\title[Large cone angles on a punctured sphere.]
{Geodesic intersections and isoxial Fuchsian groups.}

\author[G. McShane]{%
	Greg McShane 
} 
\address{%
	UFR de Math\'ematiques \\
	Institut Fourier 100 rue des maths \\
	BP 74, 38402 St Martin d'H\`eres cedex, France
}
\email{%
	Greg.McShane@ujf-grenoble.fr
}



\begin{abstract} 
The set of axes of hyperbolic elements in a Fuchsian group depends on the  commensurability class of the group.
In fact, it has been conjectured that it determines the commensurability class 
and this has been verified in for groups of the second kind by G. Mess and 
for arithemetic groups by by D. Long and A. Reid.
Here we show that the conjecture holds for almost all Fuchsian groups
and explain why our method fails for arithemetic groups.

\end{abstract}

\maketitle

\section{Introduction}
Let $\MM$ be a closed orientable hyperbolic surface. 
The free homotopy classes of closed
geodesics on $\MM$ 
conjugacy classes of hyperbolic elements in $\Gamma$.
 If  $\gamma \in  \Gamma$
 is a hyperbolic
element, 
then associated to 
$\gamma$
 is an axis 
$\axg \subset \HH$.
The projection of 
$\axg$ to $\MM$
  determines a
closed geodesic 
 whose length is 
 $\llg$.
We shall denote the set of axes of all the
hyperbolic elements in $\Gamma$
 by $\axG$.
 It's easy to check that if  $g \in \psl$ then 
 we have the relation
 \begin{equation} \label{conj}
 \ax(g\Gamma g^{-1}) = g\,\axG.
  \end{equation}

 \subsection{Isoaxial groups}
 Following Reid \cite{reid}
 we say that a pair of Fuchsian groups  $\Gamma_1$ and $\Gamma_2$
 are \textit{isoaxial} iff 
 $\ax(\Gamma_1) = \ax(\Gamma_2)$.
One obtains a a trivial example of an isoaxial pair 
by taking $ \Gamma_1$  any Fuchsian group
and $\Gamma_2 < \Gamma_1$ 
any finite  index subgroup.
This  example  can be extended 
to a more general setting as follows.
Recall that a pair of subgroups
$\Gamma_1$ and $\Gamma_2$
are \textit{commensurable}
if $\Gamma_1 \cap \Gamma_2$
is finite index in both 
$\Gamma_1$ and $\Gamma_2$.
Thus if $\Gamma_1$ and $\Gamma_2$ are commensurable then
they are isoaxial because:
$$\ax(\Gamma_1) =  \ax( \Gamma_1 \cap \Gamma_2)
= \ax(\Gamma_2),$$
It is natural to ask whether the converse is true:
 \begin{center}
 If $\Gamma_1$ and $\Gamma_2$ are isoaxial 
 then are they  commensurable?
  \end{center}
  In what follows we shall say simply
  that the group $\Gamma_1$
  is determined (up to commensurability)
  by its axes.
 We shall show that this conjecture holds
 for almost all Fuchsian groups:
  
  \begin{theorem} \label{main thm}
 For almost every point  $\rho$ in Teichmueller space 
 of a hyperbolic surface $\MM$
 the corresponding Fuchsian representation
 the fundamental group   $\Gamma$ is   determined by its axes.
  \end{theorem}
  
\subsection{Spectra}

 We define the length spectrum of $\MM$
to be the collection of lengths $\ell_\alpha$ of 
closed geodesics $\alpha \subset \MM$
\textit{counted with multiplicity}.
In fact, since $\MM$ is compact,
the multiplicity of any value in the spectrum is finite
and moreover the set of lengths is discrete.
Let $\alpha,\beta$ be primitive closed geodesics 
which meet at a point $z\in\MM$,
we denote by 
$\alpha\angle_z\beta$,
 the angle measured in the counter-clockwise direction from
 $\alpha$ to $\beta$.
 Let $\alpha,\beta$ be primitive closed geodesics 
which meet at a point $z\in\MM$,
we denote by 
$\alpha\angle_z\beta$,
 the angle measured in the counter-clockwise direction from
 $\alpha$ to $\beta$.
 Following Mondal \cite{mondal1},\cite{mondal2}
we  define an  \textit{angle spectrum}
to be the collection of  all such  angles (counted with multiplicity).

The  length spectrum has proved useful 
in studying many problems concerning 
the geometry of hyperbolic surfaces.
The angle spectrum is very different from the length spectrum:
the  set of angles is  obviously not discrete and,
as we shall see,
the there are surfaces for which 
every value has infinite multiplicity.
However, 
when considering the question of 
whether groups are isoaxial,
the angle spectrum has a distinct advantage
for it  is easy to see that:
\begin{itemize}
\item There are isoaxial groups which 
do not have the same set of  lengths,
that is, the same angle  spectrum without multiplicities.
\item If two groups are isoaxial then 
they have the same  set of angles,
that is, the same angle  spectrum without multiplicities.

\end{itemize}
Using properties of angles 
we will deduce   Theorem \ref{main thm}
from the the following lemma
inspired by a result of G. Mess (see paragraph \ref{mess} ).

\begin{lemma}\label{easy lemma}
Define the group of automorphisms of $\axG$
to be  the group of hyperbolic isometries
which preserve $\axG$.
If  $\MM$ has a 
value in its  angle spectrum
with finite multiplicity then 
$\Gamma$ is finite index 
in the group of automorphisms of $\axG$.
\end{lemma}

It remains to prove that there are such points of $\TT$,
we show in fact that they are generic:

  \begin{theorem} 
 For almost every point  $\rho \in \TT$
there is a value in the angle spectrum 
which has multiplicity exactly one.
  \end{theorem}
  
 Our method applies provided there is some value in the angle spectrum
  that has \textit{finite} multiplicity. 
  Unfortunately,  for arithemetic surfaces, the multiplicity of
  every value is infinity (Lemma \ref{infinite mult statement}).

\subsection{Sketch of proof}   The method of proof of Theorem \ref{main thm}
  follows the proof of
  the first part of Theorem 1.1 in  \cite{mcp}:
This  says that
 the set of surfaces  in Teichmeuller space 
 where every value in the  simple length spectrum has multiciplity exactly one is dense and its complement is measure zero
 ( for the natural measure on Teichmueller space.)

\subsubsection{Two properties of (simple) length functions}
Recall that   the \textit{simple length spectrum} is defined to be the collection of
lengths of simple closed  geodesics counted with multiplicity. 

There are two main ingredients used  in \cite{mcp}  :
\begin{itemize}
\item The analyticity of the geodesic length  $\ell_\alpha$
as a function over Teichmeuller space;
\item
The fact that if $\alpha,\beta$ are a pair of distinct
simple closed  geodesics
then the difference $\ell_\alpha - \ell_\beta$ 
defines a non constant 
(analytic)
 function on the  Teichmeuller space $\TT$.
\end{itemize}
It is clear that the set of of surfaces where every value in the  simple length spectrum has multiciplity exactly one
is the complement of
$$ Z := \cup _{(\alpha, \beta)}\{  \ell_\alpha -  \ell_\beta  = 0\},$$
where the union is over all pairs $\alpha, \beta$ of distinct closed simple geodesics.
Each of the sets on the left 
is nowhere dense and  its intersection with any open set  is measure zero.
Since  $Z$ is countable union of such sets,
its complement is dense and meets every  open set in a set of full measure.

We note in passing that the second of these
properties is not true without the hypothesis "simple".
Indeed, there are pairs of 
distinct closed unoriented  geodesics $\alpha\neq \beta$
such that $\ell_\alpha = \ell_\beta$ identically on $\TT$
(see \cite{Buser} for an account of their construction).

\subsubsection{Analogues for angles}

We will deduce Theorem \ref{main thm} 
using the same approach but
instead of geodesic  length functions 
we  use  angle functions.
The most delicate point is to show that if
$\alpha_1, \alpha_2$ are a pair of simple  closed  geodesics
that meet in a single point $z$
and 
$\beta_1, \beta_2$ are a pair of  closed  geodesics
that meet in a point $z'$
then the difference 
$\alpha_1\angle_z \alpha_2 -\beta_1\angle_{z'}\beta_2 $
defines a non constant function on Teichmueller space.

We do this by establishing the analogue of the following
property of geodesic length functions: 
\begin{fact}
A closed geodesic $\alpha \subset \MM$ is simple if and only if the the image of the  geodesic length function $\ell_\alpha$ is $]0,\infty[$. 
\end{fact}
Our main technical result (Theorem \ref{generic intersection})  
is an analogue of this property.
We  consider 
pairs of simple closed geodesics
 $\alpha_1, \alpha_2$
 which meet in a point $z$ --
 this configuration will be the analogue 
 of a simple closed geodesic.
 Now, for any such pair we find  a subset $X \subset \TT$
such that, for any other pair of closed geodesics
$\beta_1, \beta_2$ which meet in $z' \neq z$:
\begin{itemize}
\item the image of $X$  under  $\beta_1\angle_{z'}\beta_2$
is a proper subinterval of $]0,\pi[$ 
\item whilst its image under 
 $\alpha_1\angle_z \alpha_2$
 is the whole 
of $]0,\pi[$.
\end{itemize}

\subsection{Further remarks}
Since one objective of this work is to compare systematically the properties of geodesic length and angle functions we   include an exposition  of geodesic length functions and give an account of  the characterisation of simple geodesics mentioned above our Proposition \ref{characterise simple}.

Mondal \cite{mondal1} has obtained a rigidity result 
by using a richer collection of  data than we use here.
He defines a \textit{length angle spectrum} 
and proves that this determines a surface up to isometry.
However, the set of axes does not determine the lengths of closed geodesics
and so commensurability is the best one can hope for in the context we consider here.

In paragraph \ref{infinite mult} we answer a question of Mondal in \cite{mondal2}
concerning multiplicities by observing that arithemetic surfaces are 
very special: the multiplicity of any angle  in the angle spectrum is infinite.

 \section{Automorphisms and commensurators}
To study this question we define, following Reid,
two auxilliary groups.
The first is the 
\textit{group of automorphisms of $\axG$}:
 $$\auta 
:= \{\gamma \in\psl,\, \gamma(\axG) =\axG \}.
$$
The second is  the \textit{commensurator of $\Gamma$} 
defined as:
$$
\text{Comm}(\Gamma) :=
\{\gamma \in\psl: \gamma\Gamma \gamma^{-1}
\text{ is directly commensurable with } \Gamma
\}.
$$
We leave it to the reader to check that 
$\auta$ and $\comm(\Gamma)$ are  indeed groups 
and that they contain $\Gamma$ as  a subgroup.
In fact any element $\gamma \in \comm(\Gamma)$
is an  automorphism of   $\axG$. 
To see this,  
if $\gamma \in \comm(\Gamma)$, then
$\Gamma$ and $\gamma \Gamma \gamma^{-1}$ are commensurable  so  are isoaxial.
Now by (\ref{conj}) 
one has
$$\axG = \ax(\gamma \Gamma x^{-1} ) = \gamma\axG$$ so 
 $\gamma \in \auta$.
 In summary one has a chain of inclusions of subgroups:
 $$ \Gamma < \comm(\Gamma) < \auta < \psl.$$
We shall be concerned with two cases:
\begin{enumerate}
\item $\Gamma$ is finite index in $\auta$.
\item $\auta$ is dense in $\psl$ 
so that $\Gamma$ is necessarily an  infinite index subgroup.
\end{enumerate}
The first case arises for the class of Fuchsian groups of the second kind studied by G. Mess and the second for arithemetic groups.

\subsection{Fuchsian groups of the second kind}
\label{mess}

G. Mess in an IHES preprint studied a variety of questions
relating  to $\axG$ notably proving the following result
  
  \begin{theorem}[Mess]
 If $\Gamma_1$ and $\Gamma_2$ 
  are isoaxial Fuchsian groups of the second kind
  then    they are commensurable.
    \end{theorem}

The proof of this result 
is a consequence of the fact  that,
under the hypotheses,
$\auta$ is a discrete, convex cocompact Fuchsian group.
It is easy to deduce from this that $\Gamma$
is finite index in $\axG$.

 To show that $\auta$ is discrete 
 it suffices to find a discrete subset of $\HH$,
 containing at least two points, 
 on which it acts.
Recall that  the 
\textit{convex hull of the limit set}  of $\Gamma$,
 is a convex subset $\CH \subset \HH$.
If $\Gamma$ is a Fuchsian groups of the second kind
then its limit set $\LL$ is a proper subset of $\partial \HH$
and $\CH$ is a proper subset of $\HH$
whose frontier $\partial \CH$ consists of countably many 
complete geodesics which we call \textit{sides}.
 By definition 
  $\axG$ is $\auta$-invariant 
  and so  $\CH$ is too
since, in fact,
 it is the minimal convex set containing $\axG$. 
 Now choose a minimal length perpendicular $\lambda$
 between edges  of $\CH$;
 such a minimising perpendicular exists 
 because the double of 
 $\CH / \Gamma$ is a compact surface without boundary,
 every perpendicular between edges of $\CH$
 gives rise to a closed geodesic on the double
 and the length spectrum of the double is discrete.
 Let $L$ be the  $\auta$-orbit of $\lambda$
 and observe that  $L \cap \partial \CH$ is a discrete set
 which contains at least two points.

 \subsection{Arithemetic groups}
 
 In the case of Fuchsian groups of the first kind
 Long and Reid \cite{AR} proved the conjecture for arithemetic groups.

  \begin{theorem}[Long-Reid]\label{RL}
 If a Fuchsian group  is arithmetic then 
its commensurator is exactly the 
group of automorphisms of the group.
    \end{theorem}
 Also notice that if 
$\Gamma_1$ and $\Gamma_2$
 are isoaxial Fuchsian groups, 
 then for any 
 $\gamma \in \Gamma_2$
 $$\ax(\Gamma_1) =  \ax(\gamma\Gamma_1\gamma^{-1}),$$
   and therefore
  $ \gamma \in \auta. $
  Hence
  $ \Gamma_2 < \auta.$

So by  the above discussion
 $\Gamma_2 < \comm(\Gamma_1)$,
 and  if $\Gamma_2$ is also arithmetic, 
 then $\Gamma_1$ and $\Gamma_2$ are commensurable. 
Thus they obtain as a corollary:

  \begin{coro}
  Any pair of isoaxial arithmetic Fuchsian groups is commensurable.
  \end{coro}
    
\subsubsection{Multiplicities for arithemetic groups}
\label{infinite mult}
Let $\Gamma$ be an arithemetic Fuchsian group since its commensurator 
is dense in  $\sl$ set of geodesic intersctions is ``locally homogenous" in the 
following sense:

\begin{lemma}
Let $\theta = \alpha \angle_z \beta$ be an intersection of closed geodesics
then for any open subset  
$U \subset \MM$  
there is  a pair of closed geodesics 
$\alpha_u, \beta_u$ such that:
 $$\alpha_u \angle_{z_u} \beta_u, \, z_u \in U.$$
\end{lemma}
\proof
Choose  hyperbolic elements $a,b \in \Gamma$ 
such that the axis of $a$ (resp. $b$) is a lift of 
$\alpha$ (resp $\beta$)  to $\HH$
and so that the axes meet in a  lift $\hat{z} \in \HH$ of $z$.
Since $\comm(\Gamma)$ is dense in $\sl$,
there is some element $g \in \comm(\Gamma)$ 
so that $g(\hat{z}) \in \hat{U}$ for some lift of $U$ to $\HH$.
By the commensurability of the groups  $\Gamma$ and $g\Gamma g^{-1}$
 there is a positive integer $m$ such that
$(gag^{-1})^m,( gbg^{-1})^m \in \Gamma$ 
so that the axes of these elements project to closed geodesics
$\alpha_u, \beta_u$ on $\MM$ meeting in a point $z_u$
as required.
\hfill $\Box$

An immediate corollary is:
\begin{coro}\label{infinite mult statement}
The multiplicity of any angle $\theta$
in the spectrum of an arithemetic surface $\MM / \Gamma$
is infinite.
\end{coro}

\section{Functions on Teichmeuller space}

Recall that  the Teichmeuller 
of a surface $\MM$, $\TT$, is the
set of marked complex structures
and that, by Rieman's Uniformization Theorem,
this is identified with a component of the character variety
of $\psl$-representations of $\fund$.
Thus we think of a point $\rho \in \TT$
as an equivalence class of $\psl$-representations 
of $\fund$.
We remark  that 
$\psl := \sl /  \langle -I_2 \rangle$
so that 
although  the trace $\tr \rho(a)$ is not 
well defined 
for $a\in \fund$,  
the square of the trace  $\tr^2 \rho(a)$ is 
and  so is  $| \tr \rho(a)  |$.
In fact, there is a natural topology  $\TT$  such that
for each $a \in \fund$, 
$\rho \mapsto \tr^2 \rho(a)$ is a real analytic function.

\subsection{Geodesic length}

If $a\in \fund$ is non trivial 
then there is a unique oriented closed simple geodesic
$\alpha$ in the conjugacy class $[a]$  determined by $a$.
The length of $\alpha$,  
measured in the Riemannian metric on $\MM = \HH/\rho(\fund)$),
can be computed from 
$\tr \rho(a)$ using the well-known formula
\begin{equation}\label{length formula}
| \tr \rho(a)  |=  2 \cosh( \ell_\alpha/2).
\end{equation}
There is a natural function, 
$$\ell :\TT \times\{  \text{ homotopy classes of loops}  \} \rightarrow \,]0,+\infty[$$
 which
takes the pair $\rho,[a]$ to the length 
$\ell_\alpha$ of the
geodesic in the homotopy class $[a]$. 
It is an abuse, though common in the
literature, to refer merely to {\it the length of the geodesic
$\alpha$} (rather than, more properly, the length of the geodesic
in the appropriate homotopy class).

We define the length spectrum of $\MM$
to be the collection of lengths $\ell_\alpha$ of 
closed geodesics $\alpha \subset \MM$
\textit{counted with multiplicity}.
In fact, since $\MM$ is compact,
the multiplicity of any value in the spectrum is finite
and moreover the set of lengths is discrete.

\subsubsection{Analyticity}
A careful study of properties
of length functions was made in \cite{mcp} 
where  one of the key ingredients 
is the analyticity of this class of  functions:
\begin{fact}\label{anal length}
For each closed geodesic $\alpha$, the function 
$$\TT \rightarrow ]0,+\infty[,\, \rho \mapsto\ell_\alpha$$
 is a non constant, real analytic function.
\end{fact}
See \cite{abbook} for a proof of this.
Note that,
to prove that such a function  is non constant,
it is natural to consider two cases 
according to whether the geodesic 
 $\alpha$ is simple or not:
\begin{enumerate}
\item 
if $\alpha$ is simple then by including it as a
curve in a pants decomposition one can view
$\ell_\alpha$ as one of the Fenchel-Nielsen coordinates
so it is obviously non constant 
and, moreover, takes on any value in 
$]0,+\infty[$
\item 
if $\alpha$ is not simple then it
suffices to find a closed simple geodesic
$\beta$ such that $\alpha$ and $\beta$ meet
and use the inequality (see Buser \cite{buser}) 
\begin{equation} \label{collar}
\sinh(\ell_\alpha/2)\sinh(\ell_\beta/2) \geq 1
\end{equation}
to see that if $\ell_\beta \rightarrow 0$
then $\ell_\alpha \rightarrow \infty$
and so is non constant.
\end{enumerate}

\subsubsection{Characterization of simple geodesics}
There is always a simple closed geodesic shorter 
than any given closed geodesic.
More precisely, if $\beta \subset \MM$ is a closed geodesic
which  is not simple then  by doing surgery at
the double points  one can construct a simple closed geodesic
$\beta' \subset \MM$  with $\ell_{\beta'} < \ell_\beta$.

For $\epsilon > 0$ define the
\textit{ $\epsilon$-thin part of the  Teichmeuller space $\TT$}
to be the set 
$$ \TTE : = \{ \ell_\beta < \epsilon,\, 
\forall \beta \text{ closed simple} \}
 \subset \TT.$$
 By definition, on the complement of the thin part 
 $\ell_\beta  \geq \epsilon$   for all simple closed geodesics
and since, by the preceding remark, 
there is always a simple closed geodesic shorter 
than any given closed geodesic,
 $\ell_\beta  \geq \epsilon$
 for all  closed geodesics.

\begin{prop}\label{characterise simple}
Let $\MM$ be a finite volume  hyperbolic surface.
Then a closed geodesic $\alpha \subset \MM$ is simple if and only if the infimum over $\TT$ of 
the geodesic length function $\ell_\alpha$ is zero.
\end{prop}

\proof
In one direction,
if $\alpha$ is simple then 
$\ell_\alpha$  is one of the Fenchel-Nielsen coordinates
for some pants decomposition of $\MM$
so there is some (non convergent) sequence $\rho_n \in \TT$
such that $\ell_\alpha \rightarrow 0$.

Now suppose that $\alpha$ is not simple
and we seek  a lower bound for its length.
There are two cases depending on whether 
there exists a closed simple geodesic $\beta$ 
disjoint from $\alpha$ or not.
If there is no such geodesic then $\alpha$ meets 
every simple closed geodesic $\beta \subset \MM$
and it is cusomary to call such a curve a \textit{filling curve}.
Choose $\epsilon > 0$ and consider
the decomposition of the Teichmeuller space
into the $\epsilon$-thin part and its complement.
On the thick part $\ell_\alpha \geq  \epsilon$ 
whilst on the thin part,
 by the inequality (\ref{collar}),
it is bounded below by 
$\mathrm{arcsinh} (1/ \sinh(\epsilon/2) )$.

If there is an essential  simple closed geodesic 
disjoint from $\alpha$ then 
we cut along this curve
to obtain a possibly disconnected surface with geodesic boundary.
We repeat this process to construct a compact  surface
$C(\alpha)$ such that 
$\alpha$ is a  filling curve in 
$C(\alpha)$.
By construction $C(\alpha)$ embeds isometrically 
as a subsurface of $\MM$
and since $\alpha$ is not simple $C(\alpha)$ is not an annulus.
On the other hand, by taking the Nielsen extension of $C(\alpha)$
then  capping off with a punctured disc
we obtain a conformal  embedding
$C(\alpha) \hookrightarrow C(\alpha)^*$ 
where $C(\alpha)^*$ is a punctured surface
with a natural Poincar\'e metric.
By the Ahlfors-Pick-Schwarz Lemma 
there is a  contraction between
the  metrics induced on 
$C(\alpha)$ from the metric on $\MM$ and 
 from the Poincar\'e metric on 
$C(\alpha)^*$.
A consequence of this is that 
the  geodesic in the homotopy class 
determined by $\alpha$
on $C(\alpha)$  is longer than 
the one   in $C(\alpha)^*$.
So, to bound $\ell_\alpha$ it suffices to bound 
the length of every filling curve on a punctured surface.
There are two cases. 
\begin{itemize}
\item 
If $C(\alpha)$ has an essential simple closed curve then
we have already treated this case above.
\item
If $C(\alpha)$ has no essential simple closed curves 
then it is a  3 punctured sphere
an  the bound is trivial since the Teichmueller space 
consists of a point.
\end{itemize}

\hfill $\Box$

\section{Fenchel-Nielsen twist deformation}
\label{fn review}

Whilst make no claim as to the originality 
of the material in this section 
it is included to set up  notation
give an exposition of two results
which we use in Section \ref{variation angle}.

\subsection{The Fenchel-Nielsen twist}
We choose a simple closed curve $\alpha \subset \MM$.
Following \cite{kerk}, cut along this curve, 
and take the completion of the resulting
surface  with respect to the path metric
 to obtain a  possibly disconnected
 surface with geodesic boundary $\MM'$.

 Obviously, one can recover the original surface from $\MM'$ by identifying pairs of points of  one from each of the boundary components.
  More generally, if $t \in \RR$ then a 
  \textit{(left) Fenchel-Nielsen twist along $\alpha$}
allows one to construct a new surface $\MM_t$,
 homeomorphic to $\MM$
  by identifying  the two  boundary components 
  with a left twist of distance $t$,
  i.e. the pair of points which  are identified to obtain $\MM$
  are now  separated by distance $t$ along the image of $\alpha$ in  $\MM_t$.
 Thus this construction gives rise to a  map,
   which we will call the \textit{time $t$ twist along $\alpha$},
  $$\fntb : \MM   \rightarrow  \MM_t,$$
 discontinuous for $t\neq 0$ and mapping
  $\MM \setminus \alpha$ isometrically onto 
  $\MM_t \setminus \alpha$.
 Note that $\fntb$ is 
 not unique 
 but this will not be important for our analysis,
 what
  is important, and easy to see from the construction,
  is  that the geometry 
  of   $\MM_t \setminus \alpha$ does not vary with $t$
  as we will exploit this to obtain our main result.


\subsection{The lift of the twist to  $\HH$}
\label{lift of twist}
Let $\Gamma$ be  Fuchsian group such that 
$\MM := \HH / \Gamma$ is a closed surface,
$\alpha \subset \MM$  a non separating simple closed geodesic and
 $x \not \in \alpha$ a basepoint for $\MM$.
 Now let   $A \subset \HH $  denote the set of all lifts of $\alpha$
 and $\hat{x}\in \HH$  a lift of $x$.
Then the complement of $A$ consists of
an infinite collection of  pairwise congruent, convex sets.
Moreover,  if $P$ denotes
the connected component of the complement of $A$
containing $\hat{x}$, 
then $P$ can be identified with the universal cover of
the surface $\MM \setminus \alpha$
and the subgroup $\Gamma^P < \Gamma$ that preserves $P$
is isomorphic to the fundamental group of this subsurface.
Since the geometry of $\MM_t$ does not change with $t\in \RR$
the geometry of $P$ does not change either.
This observation is the key to establishing uniform bounds 
in  the proof of Theorem \ref{generic intersection}.

Each of the other connected components of $\HH \setminus P$ 
can be viewed  as a translate of $g_i(P)$ for  some  element 
$g_i$ of  $\Gamma$ and so $\HH$ is
\textit{ tiled} by copies of $P$.
Let us consider how this tiling evolves under the time $t$  twist 
$\fntb$ along $\alpha$.
There is a unique lift $\fnt  : \HH \rightarrow \HH$
which fixes $\hat{x}$ and hence $P$.
We can calculate the image of
a translate of  $P$ 
under the lift of $\fnt$ 
by a recursive procedure.
Suppose that for some  $g_1,\dots g_n \in \Gamma$;
\begin{itemize}
\item $\cup g_i(\bar{P})$ is connected,
\item we have determined the images of  $g_1(P),\ldots g_n(P)$.
\end{itemize}
Let $g_{n+1}(P)$ be a translate of $P$
such that  
$g_{n+1}(\bar{P}) \cap  g_n(\bar{P}) = \hat{\alpha}.$
and  we consider  two cases:
\begin{enumerate}
\item If $g_n(P) = P$ 
then the   image of $g_{n+1}(P)$ is 
$\phi^t( g_{n+1}(P))$
where $\phi^t$ is a hyperbolic translation of length $t$
with axis $ \hat{\alpha}$.
\item If  $g_n(P) \neq P$ and its  image under $\fnt$ is $h(P)$
then the   image of $g_{n+1}(P)$ is 
$h\circ\phi^t\circ  g_n^{-1} (g_{n+1}(P)))$
where $\phi^t$ is a hyperbolic translation of length $t$
with axis $ g_n^{-1}(\hat{\alpha}) \subset A$.
\end{enumerate}

This procedure allows us to prove the following:

\begin{lemma} \label{converge to limit set}
Let $\LL^P \subset \partial \HH$ denote the limit set of $\Gamma^P$.
Then $\fnt$ admits a  canonical extension 
 $\fnt : \HH \sqcup \partial \HH \rightarrow  \HH \sqcup \partial \HH$
 which is continuous on $\partial \HH$.
 Further: 
\begin{enumerate}
\item For any $w \in \Gamma^P$ one has
$\fnt(w)  = w$;
\item For any $w \in \partial \HH$ one has
$ \lim_{t \rightarrow \pm \infty} \fnt(w) \in \LL^P$
and further this is an endpoint of an edge of $\BP$.
\end{enumerate}

\end{lemma}

\proof 
It is standard from the theory of negatively curved  groups 
that the  lift admits a unique extension to $\HH \sqcup \partial \HH$,
continuous on the boundary $\dHH$,
since $\HH/ \Gamma$ is compact and
so  the restriction of  the lift to 
the set of lifts of a base point $x \in \MM$,
 $\Gamma.\{ \hat{x}\}$
is Lipschitz.

Since the extension is continuous, 
to prove (1)   it suffices to note that 
the lift of the Fenchel-Nielsen deformation  fixes the endpoints of the edges of $\BP$ and these are dense in $\LL^P$.

For (2) let  $w \in \partial \HH$ and suppose that 
it is not a point of $\BP$. 
Then there is an edge  $\hat{\alpha}$ of $\BP$
such that $w$ is a point of the interval 
determined by the endpoints of this geodesic.
It is easy to check using our recursive description 
of the action of $\fnt$ on  $\HH$
that $w$ converges to the appropriate endpoint of 
 $\hat{\alpha}$.
 \hfill $\Box$
 
 We note that (2) can also be proved as follows.
 For $t = n \ell_\alpha, n \in \ZZ$ the 
 Fenchel-Nielsen twist coincides with a Dehn twist.
If $\beta$ is a loop, disjoint from $\alpha$
then (up to homotopy)  it is fixed  by the Dehn twist.
If $\beta$ is a loop which crosses $\alpha$
then under iterated Dehn twists $\tw^n$ it limits 
to  a curve on $\MM$ that spirals to $\alpha$.
That is, lifting to $\HH$ and considering the
extension of the lift of the Dehn twist
 $\tw^n : \HH \sqcup \dHH \rightarrow \HH \sqcup \dHH $,
an endpoint of $\tw^n(\beta)$ converges to 
an endpoint of some lift of $\alpha$.
It is not difficult to pass to general $t$
using the fact that the $\fnt$ extends to a 
homeomorphism on $\HH \sqcup \dHH$.


\subsection{Separated geodesics}

We say that a pair of 
geodesics  $\hat{\gamma}_1,\hat{\gamma}_2 \subset \HH$
are \textit{separated by a a geodesic} $\hat{\gamma}$
with end points $\hat{\gamma}^\pm \in \dHH$
if the ideal  points of  $\hat{\gamma}_1,\hat{\gamma}_2$
are in different connected components of
$\dHH \setminus \{ \hat{\gamma}^\pm \}.$
Note that $\hat{\gamma}_1,\hat{\gamma}_2$
are necessarily disjoint.

If $\gamma_1, \gamma_2 \subset \HH$ are
a pair of simple closed geodesics,
such that $\alpha,\gamma_1,\gamma_2$ are disjoint 
and
 we choose an arc $\beta$ between  $\gamma_1$ and $\gamma_2$
that meets $\alpha$ transversely  in a single point
then this configuration lifts to  $\HH$
as $\hat{\gamma}_1,\hat{\gamma}_2$ separated by 
a lift $\hat{\alpha}$ of $\alpha$.
It is easy to convince oneself that,
 as we deform by the Dehn twist $\tw^n$,
 the length of $\beta$ goes to infinity.
Essentially, our next lemma says that this is true 
 for any pair of geodesics $\gamma_1, \gamma_2$
  in $\MM$ admitting an arc that meets $\alpha$
  in an essential way.

\begin{lemma}
Let $\hat{\gamma}_1,\hat{\gamma}_2 \subset \HH$ be a pair of geodesics which are separated by some lift of $\alpha$
then the distance between 
$\fnt(\hat{\gamma}_1)$ and $\fnt(\hat{\gamma}_2)$ tends to infinity 
as $t\rightarrow \pm \infty$.
\end{lemma}

\proof
Let $\hat{\alpha}$ be a lift of $\alpha$ which 
separates $\hat{\gamma}_1,\hat{\gamma}_2 \subset \HH$.
Let $P_1$ and $P_2$ be the pair of 
complementary regions which have $\hat{\alpha}$
as a common edge
and we label these so that
$\hat{\gamma}_i$ is on the same side of  $\hat{\alpha}$
as $P_i$ for $i=1,2$.
We choose the lift of the base point to be in $P_1$
and lift the Fenchel-Nielsen deformation.

First consider the  orbit $\fnt(y)$ of an ideal endpoint $y$ 
of  $\hat{\gamma}_2$ as $t \rightarrow \infty$.
Since $x \in P_1$, the region $P_2$ 
gets translated  and so, for any side  $\beta$ of $P_2$,
the sequence $\fnt(\beta)$ converges to 
the endpoint $\hat{\alpha}^+$.
Now there is a pair of edges $\beta_1, \beta_2$
such that the endpoints of $\hat{\gamma}_2$
are contained in the closed interval 
containing the endpoints of $\beta_1, \beta_2$.
Since each of the $\beta_i$ converge to $\hat{\alpha}^+$
under the deformation it is  easy to see that $\fnt(\hat{\gamma}_2)$
must converge to $\hat{\alpha}^+$ too.

Now consider the orbit of an endpoint $y$ 
of  $\hat{\gamma}_1$ 
under the deformation.
It suffices to show that,
under this deformation, 
$y$ does not converge 
to $\hat{\alpha}^+$.
There are two cases
according to whether or not $y$ belongs to the limit set $\LL^{P_1}$
of the  subgroup of $\Gamma$  which stabilises $P_1$.
\begin{enumerate}

\item If  $y \in \LL^{P_1}$ then it is invariant under the 
Fenchel-Nielsen deformation.
\item If  $y \not \in \LL^{P_1}$ then it limits to a 
point in $y_\infty \in \LL^{P_1}$ which is an endpoint of one of
the edges of $P_1$. 
By hypothesis 
$\hat{\gamma}_1$ does not meet $\hat{\alpha}$ and
so $y_\infty$ is not $\hat{\alpha}^+$.

\end{enumerate}

\hfill $\Box$

%

\section{Geodesic angle functions}

We present two methods for computing 
(functions of) the angle $\aza$ between 
$\alpha_1,\alpha_2$ at $z$.
The first method,
just like the formula (\ref{length formula}) for geodesic length,
is a closed  formula 
in  terms of traces  (equation (\ref{relation})
whilst the second
 is in  terms of end points of 
lifts of $\alpha_1,\alpha_2$ to the Poincar\'e disk
 (equation (\ref{CR formula})).
This second formula will prove useful
for obtaining estimates for the variation
of angles along a Fenchel-Nielsen deformation.
In either case, 
 we start as befor by  identifying $\MM$ with 
the quotient $\HH/ \Gamma$ where 
$\Gamma = \rho(\fund),\, \rho \in \TT$.
We choose $z$ as a basepoint for $\MM$
and associate elements $a_1,a_2 \in  \pi_1(\MM,z)$
such that $\alpha_i$ is the unique oriented closed geodesic
in the conjugacy class $[a_i]$
in the obvious way.

\subsection{Traces and analyticity}
As explained in the introduction we shall need an analogue of Fact \ref{anal length}
 so we give a brief account of the analyticity of the angle functions:

\begin{prop}
If  $\rho \in \TT$
is a point in Teichmuller space then 
$$ \TT \rightarrow ]0, 2\pi[,\,
\rho \mapsto \aza,$$
is a  real analytic function.
\end{prop}

\proof 
With the notation above
 we have the following expression for the angle:
\begin{equation}\label{relation}
 \sin^2(\alpha_1 \angle_z \alpha_2) =   
 \frac{
 4( 2 -  \mathrm{tr}
  [\rho(a_1),\rho(a_2) ])
  }
  {( \mathrm{tr}^2 \rho(a_1)-4)( \mathrm{tr}^2 \rho(a_2)-4)}.
 \end{equation}
This equation is actually implicit in  \cite{pjm}
 but it is not claimed to be new there and seems to have been well known.
 The left hand side of  (\ref{relation})
 is clearly an analytic function on $\TT$
 and it follows from elementary real analysis
 the the angle varies real analytically too.
 \hfill $\Box$
 
Note that, though we will not need this, 
(\ref{relation}) shows that 
the square of the sine  is in fact a rational function of traces (see Mondal \cite{mondal2}
for applications of this).

 \subsubsection{Cross ratio formula}
 It will useful to to have another formula for the angle
 in terms of a cross ratio .
 This formula is well-known, see for example,
\textit{ The Geometry of Discrete Groups},
 by A.F. Beardon but
 we since we will use it extensively to obtain bounds
 we give a short exposition.
If $\theta$ is the angle between two hyperbolic geodesics
$\hat{\alpha},\hat{\beta} \subset \HH$
then $\tan^2 (\theta/2)$ can be expressed as a cross ratio.
One can prove this directly by taking
 $\hat{\alpha}$ to 
have endpoints 
$\alpha^\pm = \pm 1$ and $\hat{\beta}$ endpoints 
$\beta^\pm = \pm e^{i\theta}$ 
in the Poincar\'e disc model.
Then
\begin{equation}\label{CR formula}
 \left( \frac{\alpha^+  -  \beta^+  }{\alpha^+  -  \beta^-  } \right)
 \left( \frac{\alpha^-  -  \beta^ -  }{\alpha^-  -  \beta^+ } \right)
 = 
\left( \frac{1 - \ei }{1  + \ei} \right)
\left( \frac{ -1  + \ei}{ -1 - \ei } \right)
=  \left( \frac{1 - \ei }{1  + \ei} \right)^2 = \tan^2 (\theta/2).
\end{equation}


\section{Angles defined by closed geodesics}


\subsection{Variation of angles}
\label{variation angle}

In this paragraph we give an improved version
of the following well known fact:

\begin{fact}\label{full range1}
Let  $\alpha,\beta \subset \MM$
be a pair of   closed simple geodesics
that meet in a point
$z\in \MM$.
If $\alpha$ is simple then  for any $\theta \in ] 0,\pi[$ 
there exists $\rho \in \TT$ such that
$$\alpha\angle_x\beta = \theta.$$ 
\end{fact}
Under the hypothesis,  there is a convex
subsurface $\MM' \subset \MM$
homeomorphic to a holed torus
which contains $\alpha \cup \beta$.
The fact follows   by presenting
$\MM'$ as the quotient of $\HH$ by a Schottky group.

Using the preceding
discussion of the Fenchel-Nielsen deformation
we can relax the hypothesis on $\beta$
even whilst taking the restriction of the angle function 
to a one dimensional submanifold of $\TT$.
The proof will should also 
serve to familiarise the reader with the notation
and provide intuition as to why this case
is different to that of an intersection 
of a generic pair of closed geodesics
treated in Theorem \ref{generic intersection}

\begin{lemma}\label{full range}
Let  $\alpha,\beta \subset \MM$
be a pair of   closed geodesics
that meet in a point
$z\in \MM$.
If $\alpha$ is simple then  for any $\theta \in ] 0,\pi[$ 
and any  $\rho_0 \in \TT$ there exists
$\rho_t \in \TT$ obtained from $\rho_0 $ 
by a time $t$ Fenchel-Nielsen twist along $\alpha$
such that 
$$\alpha\angle_x\beta = \theta.$$ 
Moreover,  
$$ \lim_{ t \rightarrow \pm \infty}
 \alpha\angle_x\beta \in \{ 0, \pi\}. $$
\end{lemma}

\proof

With the notation of subsection \ref{lift of twist},
there  is a convex region $P$ in $\HH$ bounded 
by lifts of $\alpha$ as before.
Let $\hat{\alpha}$   be an edge of $\BP$, 
and choose a corresponding lift $\hat{\beta}$ which intersects
$\hat{\alpha}$.
There is an element of the covering group $g \in \Gamma$ such that 
$$\hat{\alpha} = \BP \cap g(\BP).$$
We lift the Fenchel-Nielsen deformation 
and consider, as before, its extension
$$\fnt : \HH \sqcup  \partial \HH \rightarrow \HH \sqcup  \partial \HH .$$
Now, arguing as in Lemma \ref{converge to limit set},
we see that;
\begin{itemize}
\item 
the endpoints of $\hat{\alpha}$ are fixed by $\fnt$,
\item the endpoint of 
$\hat{\beta}$ on the same side of $\hat{\alpha}$ as $P$
converges to a point $z \neq \alpha^+$  as $t \rightarrow -\infty$,
\item the other endpoint of
$\hat{\beta}$ converges to  $ \alpha^+$  as $t \rightarrow -\infty$
\end{itemize}
It follows that, 
after possibly changing the orientation of $\beta$, 
that the angle between
$\hat{\alpha}$ and $\hat{\beta}$,
and hence $\alpha \angle_x \beta$,
tends to $0$.
Likewise, as $t \rightarrow +\infty$
the angle between
$\hat{\alpha}$ and $\hat{\beta}$,
and hence $\alpha \angle_x \beta$,
tends to $\pi$.

Thus, by continuity, the range of the angle function is $]0,\pi[$.

\hfill $\Box$


\begin{theorem}\label{generic intersection}
Let $\beta_1,\beta_2$ a pair of closed geodesics and
 $y \in\beta_1 \cap \beta_2 $.
Then for any simple closed  geodesic $\alpha$,
different from both  $\beta_1$ and $\beta_2$,
the angle function $\beta_1\angle_y\beta_2$
is bounded away from  $\pi$
along the Fenchel-Nielsen orbit of  $\rho \in \TT$.
\end{theorem}

\proof

If $\alpha$ and 
$\beta_1\cup \beta_2$ are disjoint
then $\beta_1\angle_y\beta_2$
is constant along the  $\fnt$-orbit
so the result is trivial.

Suppose now that $\alpha$ and 
$\beta_1\cup \beta_2$ are not disjoint
and choose $x$ as a basepoint of $\MM$.
Then, with the notation of paragraph \ref{??},
there  is a convex region $P$ in $\HH$ bounded 
by lifts of $\alpha$.
 We now  consider 
three cases 
according to the number of 
edges of $\BP$ that $\hat{\beta_1} \cup \hat{\beta_2}$ 
meets.

\vspace{.2in}
\begin{figure}[h]
\label{4lifts}
\begin{center}
 \includegraphics[scale=.5]{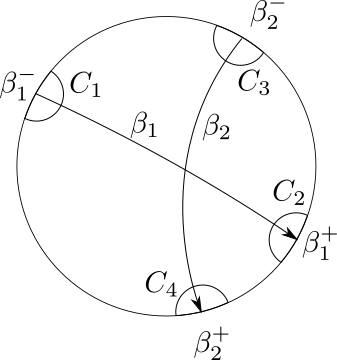} 
\end{center}
\caption{Case of 4 intersections.}

\end{figure}

We first deal with the simplest case.
Suppose that  $\hat{\beta_1} \cup \hat{\beta_2}$ 
meets  $\BP$ in four distinct edges
denoted 
$C_1,C_2,C_3,C_4 \subset \HH$,
and, after possibly relabelling these,
$\hat{\beta_1}$ meets $C_1,C_2$
whilst
$\hat{\beta_2}$ meets $C_3,C_4$
as in Figure \ref{4lifts}.
Now we deform $\rho_0$ by 
a Fenchel-Nielsen twist along $\alpha$
to obtain a 1-parameter family of
$\rho_t \in \TT, t\in \RR$.
As we have seen above, under
 such a deformation 
the length of  $\alpha$
does not change nor does
the geometry of $\BP$
in particular the positions of the  $C_i$ remain unchanged.
From our discussion of the $\fnt$ and its extension to 
 $\HH \sqcup \dHH$
it is clear that,
$\forall t \in \RR$,
$\hat{\fnt(\beta_1)}$ meets $C_1,C_2$
whilst
$\hat{\fnt(\beta_2)}$ meets $C_3,C_4$.
Thus, if the diameters of the circles were small,
 the angle at $\hat{z}$ 
cannot not vary much from its value at $\rho_0$
since the radii of the circles are small.
More generally,  we can bound the size of the angle using the cross ratio formula.
Labeling the endpoints as in Figure \ref{4lifts}
one has:
\begin{equation*}
\tan^2 (\theta/2) =
  \left| \frac{\beta_1^+  -  \beta_2^+  }{\beta_1^+  -  \beta_2^-  } \right|
 \left| \frac{\beta_1^-  -  \beta_2^ -  }{\beta_1^-  -  \beta_2^+ } \right|
\end{equation*}
Note first that each of the four points lies on the unit circle 
and so  that its  diameter, that is 2,
is a trivial upper bound for each of the four distances
appearing on the left hand side of this equation.
Now under the deformation each of the endpoints 
$\fnt(\beta_i^\pm)$ stays in one of four disjoint 
euclidean discs defined by one of the $C_j$.
In particular,  there exists  $\delta_4 > 0 $ such that
for all $t \in \RR$
\begin{eqnarray*}
\delta_4 \leq |\fnt(\beta_1^\pm)   -  \fnt(\beta_2^\pm) | & \leq & 2 \\
\delta_4 \leq  |\fnt(\beta_1^\pm)   -  \fnt(\beta_2^\mp) | & \leq &  2
\end{eqnarray*}
and this is sufficient to obtain bounds  on  the cross ratio:
\begin{equation}
1/2 \delta_4 \leq   \tan (\theta/2) \leq  2/\delta_4
\end{equation}

If $\hat{\beta_1} \cup \hat{\beta_2}$ 
meets  $\BP$ in just two  edges,
$C_1,C_2  \subset \HH$ say.
Although
 we no longer have a uniform  lower bound for  
$|\fnt(\beta_1^\pm)   -  \fnt(\beta_2^\pm) |$
 in this case,
there still  exists $\delta_2 > 0$ such that  for all $t \in \RR$,
$$\delta_2 \leq |\fnt(\beta_1^\pm)   -  \fnt(\beta_2^\mp) |.$$
Thus, for all $t \in \RR$,
\begin{equation}
   0 \leq \tan (\theta/2) \leq  2/\delta_2.
\end{equation}

\vspace{.2in}
\begin{figure}[h]
\label{2 edges}
\begin{center}
\includegraphics[scale=.5]{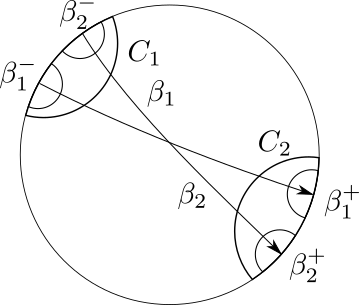} 
\end{center}
\caption{Case of 2 intersections.}
\end{figure}
Finally, if $\hat{\beta_1} \cup \hat{\beta_2}$ 
meets  $\partial P$ in exactly  three  edges
then it is easy to see that,
using the same reasoning as for the two edge case,
there is $\delta_3$ such that
$$\delta_3 \leq |\fnt(\beta_1^\pm)   -  \fnt(\beta_2^\mp) |.$$

\hfill $\Box$

\begin{coro}
Let 
 $\alpha_1, \alpha_2$
 pairs of simple closed geodesics
 which meet in a single  point $z$ 
and  
$\beta_1, \beta_2$ primitive closed geodesics which meet in $z'$.
If the difference
$$\alpha_1 \angle_z \alpha_2 - \beta_1 \angle_{z'} \beta_2$$
is constant then the angles are equal and,
after possibly relabelling the geodesics,
$\alpha_i = \beta_i$ and $z = z'$.
\end{coro}

Note that we cannot suppose that $z,z'$ are distinct
because of the following phenomenon.
Jorgenson studied intersections of closed geodesics
proving in particular that if $z \in \MM$ was the intersection 
of a pair of distinct closed geodesics then it is the intersection
of infinitely many pairs of distinct closed geodesics.
Such intersections are \textit{stable} 
in that, if $(\alpha_i)_i$ is a family of geodesics
obtained from Jorgenson's procedure 
that meet in a point $z$ on some hyperbolic surface $\MM$
then, for any $\rho \in \TT$,
 there is   $z_\rho \in \HH/ \rho(\fund)$
 such that the family  $(\alpha_i)_i$  meet in  $z_\rho$.

\proof

We first consider the case where four geodesics are distinct then,
under the Fenchel Nielsen twist  along  $\alpha_1$, 
the image of $\alpha_1 \angle_z \alpha_2$  is $]0,\pi[$
whilst, by Theorem \ref{generic intersection},
the image  $\beta_1 \angle_{z'} \beta_2$  is 
a strict subinterval. 
It is easy to see
$\alpha_1 \angle_z \alpha_2 - \beta_1 \angle_{z'} \beta_2$
cannot be a constant.
 
Now suppose, 
$\alpha_1 = \beta_1$.
if $\alpha_2= \beta_2$ then,
 since $\alpha_1$ and $\alpha_2$ meet in a single point,
 we must have $z = z'$
 and the angles must be the same.
 
On the other hand,  
if $\alpha_2 \neq  \beta_2$ 
then $z,z'$ may or may not be distinct
\begin{itemize}

\item If $z = z'$ then,
by  Lemma \ref{full range}, 
 both
$\alpha_1 \angle_z \alpha_2$ and $\beta_1 \angle_{z'} \beta_2$
tend to $0$ or $\pi$ as the Fenchel-Nielsen parameter $t \rightarrow \pm \infty$.
Therefore, if the difference is constant it must be $0$ or $\pi$
and so, up to switching orientation,
 $\beta_1 = \beta_2$.
 \item 
 If $z = z'$ then,
by   Lemma   \ref{full range} and Theorem \ref{generic intersection}
$\beta_1 \angle_{z'} \beta_2$
is a proper subinterval  of the range of 
$\alpha_1 \angle_z \alpha_2$ 
so the difference cannot be constant. 
\end{itemize}

\hfill $\Box$

\proof of Lemma \ref{easy lemma}
Suppose that $\MM$  has a value
 in its angle spectrum, $\theta$ say,  with finite multiplicity.
Let $x_1,x_2 \ldots x_n \in \MM$ be a complete list of points
such that there are pair of slosed geodesics meeting 
at $x_i$ at angle $\theta$.
Then the set of preimages    of the  $x_i$  
under the covering map $\HH \rightarrow \MM$
is a discrete set which is invariant under $\auta$.
Thus $\auta$ is discrete and 
has $\Gamma$ as a finite 
index subgroup.

\hfill $\Box$

\end{document}